\documentclass[11pt]{article}
\usepackage{epigamath}

\usepackage[notext]{kpfonts}
\usepackage{baskervald}

\setpapertype{A4}

\usepackage[english]{babel}

\usepackage[dvips]{graphicx}     

\title{Hyper-K\"ahler Fourfolds Fibered by Elliptic Products}
\titlemark{Hyper-K\"ahler Fourfolds Fibered by Elliptic Products}
\author{\vspace{0cm} Ljudmila Kamenova}
\authoraddresses{
\authordata{Ljudmila Kamenova}{\firstname{Ljudmila} \lastname{Kamenova}\\
\institution{Department of Mathematics, Room 3-115, Stony Brook University, 
Stony Brook, NY 11794-3651, USA}\\
\email{kamenova@math.stonybrook.edu}}
}
\authormark{L. Kamenova}
\date{\vspace{-5ex}} 
\journal{\'Epijournal de G\'eom\'etrie Alg\'ebrique} 
\acceptation{Received by the Editors on October 10, 2017, and in final form
on May 8, 2018. \\ Accepted on August 23, 2018.}

\newcommand{\articlereference}{Volume 2 (2018), Article Nr. 7}

\acknowledgement{Research partially supported by NSF 
Grants DMS-0111298 and DMS-0805594}

 \usepackage[all]{xy}

\allowdisplaybreaks

\newtheorem{theorem}{Theorem}[section]

\newtheorem{thm}[theorem]{Theorem}

\newtheorem{defn}[theorem]{Definition}

\begin{document}


\maketitle



\begin{prelims}


\def\abstractname{Abstract}
\abstract{Every fibration of a projective hyper-K\"ahler fourfold has fibers which 
are Abelian surfaces. In case the Abelian surface is a Jacobian of a genus 
two curve, these have been classified by Markushevich. We study 
those cases where the Abelian surface is a product of two elliptic curves, 
under some mild genericity hypotheses.}



\vspace{0.15cm}

\languagesection{Fran\c{c}ais}{%

\textbf{Titre. Vari\'et\'es hyper-k\"ahleriennes de dimension 4 fibr\'ees par des produits elliptiques} \commentskip \textbf{R\'esum\'e.} Toute fibration d'une vari\'et\'e hyper-k\"ahlerienne projective de dimension 4 a pour fibres des surfaces ab\'eliennes. 
Dans le cas o\`u la surface ab\'elienne est la jacobienne d'une courbe de genre 2, 
elles ont \'et\'e classifi\'ees par Markushevich. Nous \'etudions le cas o\`u la surface ab\'elienne est un produit de deux courbes elliptiques, sous des hypoth\`eses de g\'en\'ericit\'e peu restrictives.}

\end{prelims}


\newpage

\setcounter{tocdepth}{1} \tableofcontents

\section{Introduction} 

Among projective complex manifolds with zero first Chern class, hyper-K\"ahler 
manifolds play a distinguished role, generalizing the class of K3
surfaces.  In fact the hyper-K\"ahler manifolds of complex dimension
two are precisely the K3 surfaces.
In higher dimensions, dimension four and higher, few hyper-K\"ahler 
manifolds are known: two infinite classes introduced by Beauville 
(\cite{Beauville1}) and two exceptional cases discovered by O'Grady 
(\cite{og1} and \cite{og2}).  Nonetheless, even in
the first case of hyper-K\"ahler fourfolds, we are still far from a
classification or even a proof that hyper-K\"ahler fourfolds form a
bounded family.

Just as elliptically fibered K3 surfaces are more amenable to study, 
also fibered hyper-K\"ahler manifolds over smooth bases are better 
understood. By works of Matsushita and Hwang, (\cite{m, m2} and \cite{hw}), 
every holomorphic fibration whose total space is a $2n$-dimensional 
hyper-K\"ahler manifold has base manifold $\mathbb{P}^n$ and has general 
fiber an $n$-dimensional Abelian variety.  
A hyper-K\"ahler manifold fibered 
over $\mathbb{P}^n$ with general fiber an Abelian variety is called an 
\emph{Abelian  fibration}.  Every Abelian fibered manifold is birational to a 
``Tate-Shafarevich twist'' of an Abelian fibration with a section, and the 
set of these twists coming from a given Abelian fibration 
is classified by a well-studied Tate-Shafarevich group.  Thus, from now
on we assume that the Abelian fibration admits a section. 

The principal case is when the fibers are \emph{principally polarized}
Abelian varieties.  For hyper-K\"ahler fourfolds, the fiber will be a
principally polarized Abelian surface.  There are two types of such
surfaces, depending on whether or not the polarizing divisor is
irreducible or reducible, i.e., whether the Abelian surface is the
Jacobian of a genus $2$ curve or a product of elliptic curves.
Markushevich classified the first case, when the Abelian surface is
the Jacobian of a genus-$2$ curve.

\begin{thm}[Markushevich, \cite{markushevich96}]
Every Abelian fibered irreducible holomorphic symplectic fourfold $X$ which 
is the compactified relative Jacobian of a family of genus-two curves 
is birational to ${\it Hilb}^2 (S)$, for a K3 surface $S$ which is a 
double cover of $\mathbb P^2$ branched over a plane sextic. 
\end{thm}

Our goal is to classify those Abelian fibrations where the general 
fiber is a product of two elliptic curves, an \emph{elliptic product}.  
These hyper-K\"ahler 
fourfolds can become ``wild'' when the fibration is allowed to become
pathologically degenerate.  In fact, we conjecture that these
pathological degenerations never occur, at least for a general member
of a deformation family.  But, at the moment, we add our conjectured 
genericity conditions as hypotheses. Here is the main result of the paper. 

\begin{thm}
Consider a double cover map $f : \bar V \rightarrow \mathbb P^2$, where 
$\bar V$ is a normal variety, and a rational elliptic fibration 
$\pi : {\mathcal E} \dashrightarrow \bar V$ with a section $\sigma$, 
where ${\mathcal E}$ is projective. Let $D$ be the branch locus 
of $f$, $\tilde D = f^{-1}(D)$ and let $G$ be the discriminant locus
of $\pi$. 
Let $i$ be the involution acting on $\bar V$ that interchanges the sheets of 
the double cover $f$ and let $i_{\mathcal E}$ be the induced involution 
on ${\mathcal E} \times_{\bar V} i^\ast {\mathcal E}$.  
Assume the following genericity conditions are satisfied: 
\begin{enumerate}
\item[\rm (1)] Over the general point of $G$, the fiber of $\pi$ is
      irreducible and semistable. 
\item[\rm (2)] The image $f(G)$ contains no singular point of $D$, and $G \cap i(G)$ 
      is finite.  
\end{enumerate}
Assume that $({\mathcal E} \times_{\bar V} i^\ast {\mathcal E}) / 
i_{\mathcal E} \dashrightarrow \mathbb P^2$ admits a resolution of singularities 
$p: X \rightarrow \mathbb P^2$ which is an Abelian fibration on a 
projective hyper-K\"ahler manifold $X$ with a section $\tau$. Then 
$\bar V \cong \mathbb P^1 \times \mathbb P^1$ and $X$ is birational to 
$\text{Hilb}^2(S)$ for a projective elliptic K3 surface $S$ with a section. 
\end{thm}




Notice that for every projective, elliptic K3 surface $S$ with a
section, $\text{Hilb}^2(S)$ is a projective hyper-K\"ahler fourfold fibered
by elliptic products over $\text{Hilb}^2(\mathbb P^1) = \mathbb P^2$, and 
admits a section. In this example $\bar V = \mathbb P^1 \times \mathbb P^1 
\rightarrow \mathbb P^2$ is the standard double covering 
and the divisor ${\mathcal E} \subset \text{Hilb}^2(S) \times_{\mathbb P^2} 
\bar V$ appears naturally. The question that we study inverts this 
construction under the assumptions $(1)-(2)$ and we show that the only 
examples one obtains are birational to $\text{Hilb}^2(S)$. 
There are countably many families of projective, elliptic K3 surfaces
$S$, hence countably many families of elliptic product hyper-K\"ahler
fourfolds. 
We conjecture that 
every sufficiently general elliptic product hyper-K\"ahler fourfold 
satisfies the genericity conditions.  All such are deformation 
equivalent to one of the countably many families arising from 
projective, elliptic K3 surfaces.

\section{Preliminaries}

First we define our main objects of study, {\it irreducible holomorphic 
symplectic manifolds} or {\it irreducible hyper-K\"ahler manifolds}. 

\begin{defn}{\rm
A compact complex K\"ahler manifold $X$ is called {\it irreducible 
holomorphic symplectic}  
if it is simply connected and if $H^0(X, \Omega_X^2)$ is spanned by an 
everywhere non-degenerate 2-form  $\sigma$.}
\end{defn}

The two-form $\sigma$ is {\it everywhere non-degenerate} if it induces 
an isomorphism ${\mathcal T}_X \to \Omega_X$. 
The last condition in the definition implies that 
$h^{2,0}(X)=h^{0,2}(X)=1$ and $K_X \cong {\mathcal O}_X$, i.e., 
$c_1(X)=0$.

\begin{defn}{\rm
A compact connected $4n$-dimensional Riemannian manifold $(M,g)$ is called
{\it irreducible hyper-K\"ahler} if its holonomy group is ${\rm Sp}(n)$.}
\end{defn}

As Huybrechts notes in \cite{h}, irreducible holomorphic symplectic manifolds 
with a fixed K\"ahler class are equivalent to 
irreducible hyper-K\"ahler manifolds with a fixed complex structure. 

\begin{defn}{\rm
An {\it Abelian fibration} on a $2n$-dimensional hyper-K\"ahler manifold $X$ 
is the structure of a fibration over ${\mathbb P}^n$ whose general fiber is 
a smooth abelian variety of dimension $n$.}
\end{defn}

This is a higher dimensional analogue of elliptic
fibrations on K3 surfaces. Any fibration structure of a projective 
hyper-K\"ahler manifold is an Abelian fibration due to the following 
theorems by Matsushita~\cite{m, m2} and Hwang~\cite{hw}:

\begin{thm} [Matsushita, \cite{m, m2}]
\label{matsushita} 
For a projective irreducible holomorphic symplectic manifold $X$, let 
$f : X \rightarrow B$ be a proper surjective 
morphism such that the general fiber $F$ is connected. Assume that 
$0<{\mathrm{dim}}B<{\mathrm{dim}}X$. Then 
\begin{enumerate}
\item[\rm (1)] $F$ is an abelian variety;  
\item[\rm (2)] $B$ is $n$-dimensional. If $B$ is smooth, then it has the same Hodge 
numbers as ${\mathbb P}^n$;  
\item[\rm (3)] the fibration is Lagrangian with respect to the holomorphic
symplectic form.
\end{enumerate}
\end{thm}

\begin{thm} [Hwang, \cite{hw}]
In\,the setting of Matsushita's\,theorem, if $B$ is smooth, then $B$ is 
biholomorphic to~${\mathbb P}^n$. 
\end{thm}


Fogarty gives the following description of the Hilbert schemes of  
complex surfaces. 

\begin{thm} [Fogarty, \cite{f}] \label{fog} 
For a nonsingular surface $X$ and $n \in \mathbb N$
\begin{enumerate}
\item[\rm (1)]  ${\it Hilb}^n (X)$ is non-singular of dimension $2n$, and
\item[\rm (2)] $\pi : {\it Hilb}^n (X) \rightarrow S^n X$ is a resolution of 
singularities, where $S^n X$ is the n-th symmetric product of $X$.
\end{enumerate}
\end{thm}

Due to Beauville (\cite{Beauville1}), Hilbert schemes of K3 surfaces give 
one of the standard series of examples of irreducible holomorphic symplectic 
manifolds. The other standard series of examples, generalized Kummer varieties 
(also due to Beauville), arise from Hilbert schemes of abelian surfaces. 
In \cite{W1} Wieneck proved that if $X$ is a hyper-K\"ahler manifold 
that admits a Lagrangian fibration, the polarization type 
of the fibration is a deformation invariant of the fibration. 
As Wieneck shows, for Lagrangian fibrations the ``polarization type'' of the 
fiber is independent from the choice of an integral K\"ahler form as long 
as it is obtained as the restriction from a K\"ahler class on the ambient 
hyperk\"ahler manifold. Moreover, 
if $X$ is of $K3^{[n]}$-type, then the polarization type is principal. 
If $X$ is of generalized Kummer type, then Wieneck showed that the 
polarization type of a Lagrangian fibration depends on the connected 
component of the moduli space, \cite{W2}. 


\section{Proof of the Main Theorem}

In this section we prove our main theorem. 

\begin{thm}
Consider a double cover map $f : \bar V \rightarrow \mathbb P^2$, where 
$\bar V$ is a normal variety, and a rational elliptic fibration 
$\pi : {\mathcal E} \dashrightarrow \bar V$ with a section $\sigma$, 
where ${\mathcal E}$ is projective. Let $D$ be the branch locus 
of $f$, $\tilde D = f^{-1}(D)$ and let $G$ be the discriminant locus
of $\pi$. 
Let $i$ be the involution acting on $\bar V$ that interchanges the sheets of 
the double cover $f$ and let $i_{\mathcal E}$ be the induced involution 
on ${\mathcal E} \times_{\bar V} i^\ast {\mathcal E}$.  
Assume the following genericity conditions are satisfied: 
\begin{enumerate}
\item[\rm (1)] Over the general point of $G$, the fiber of $\pi$ is
      irreducible and semistable. 
\item[\rm (2)] The image $f(G)$ contains no singular point of $D$, and $G \cap i(G)$ 
      is finite.  
\end{enumerate}
Assume that $({\mathcal E} \times_{\bar V} i^\ast {\mathcal E}) / 
i_{\mathcal E} \dashrightarrow \mathbb P^2$ admits a resolution of singularities 
$p: X \rightarrow \mathbb P^2$ which is an Abelian fibration on a 
projective hyper-K\"ahler manifold $X$ with a section $\tau$. Then 
$\bar V \cong \mathbb P^1 \times \mathbb P^1$ and $X$ is birational to 
$\text{Hilb}^2(S)$ for a projective elliptic K3 surface $S$ with a section. 
\end{thm}



\begin{proof}

Notice that the intersection $G \cap \tilde D$ consists of finitely many 
points because of the genericity conditions. Indeed, if $G \cap \tilde D$ 
wasn't finite, then $f(G)$ and $D$ would have a whole component in common 
and the fibers of $p$ above this component would be very degenerate and not 
semistable. More precisely, if $D$ is reducible, and $f(G)$ contains one of 
its irreducible components, then $f(G)$ would contain singular points of $D$, 
namely the points of intersection of the given irreducible component with 
other components, and this is a contradiction. On the other hand, if $D$ is 
irreducible, then $f(G)$ would contain $D$, i.e., $\tilde D \subset G$. 
Since $\tilde D$ is fixed by the involution $i$, we would get a contradiction 
with the condition that $G \cap i(G)$ is finite. 
This is the place where we use the ``genericity conditions'' $(2)$.

\medskip

Let $L \subset \mathbb P^2$ be a line intersecting $D$ transversally at 
general points of $D$ and avoiding the finite set $f(G \cap \tilde D) = 
f(G) \cap D$. Denote the pre-image of $L$ in $\bar V$ by $L_V$. 
The restricted map $f_V : L_V \rightarrow L$ is 2:1. Since we have choosen 
$L$ to be general so that $L_V$ doesn't intersect $\tilde D \cap G$, 
then every fiber of 
${\mathcal E}_{L_V} \times_{L_V} i^\ast {\mathcal E}_{L_V} \rightarrow L_V$
over a point of $L_V \cap \tilde D$ is smooth. 
After pulling back our construction to $L_V$, there is a rational map:
$$L_V \times_{\mathbb P^2} X \dashrightarrow 
{\mathcal E}_{L_V} \times_{L_V} i^\ast {\mathcal E}_{L_V}$$
which by Weil's extension theorem \cite[Theorem 4.4/1]{neron} 
is regular on the smooth locus 
of $L_V \times_{\mathbb P^2} X \rightarrow L_V$ since the fibers are 
abelian varieties.

\medskip

Every section of $p : X \rightarrow \mathbb P^2$ is contained in the 
smooth locus of $ X \rightarrow \mathbb P^2$, so it pulls back to a 
curve in the smooth locus of 
$L_V \times_{\mathbb P^2} X \rightarrow L_V$ via the section 
$id_{L_V} \times_L  \tau|_L$.

\medskip

Since ${\mathcal K}_X = {\mathcal O}_X$, the relative dualizing sheaf of 
$p: X \rightarrow \mathbb{P}^2$ is 
$\omega_{X / \mathbb P^2} = p^* {\mathcal O}_{\mathbb P^2}(3)$, and 
therefore $\tau^\ast \omega_{X / \mathbb P^2} = {\mathcal O}_{\mathbb P^2}(3)$. 
The section $\sigma$ induces a section ($\sigma$, $\sigma \circ i$) of the map:
$$pr_{\bar V} : {\mathcal E} \times_{\bar V} i^\ast {\mathcal E} \dashrightarrow 
\bar V.$$

Let $\tilde \tau: \bar V \rightarrow \bar V \times_{\mathbb P^2} X$ 
be the section of $\bar V \times_{\mathbb P^2} X \rightarrow \bar V$ 
obtained by base change from $\tau$. We consider $\tau f_V (L_V)$ as a 
section of the map $L_V \times_{\mathbb P^2} X \rightarrow L_V$ and in the 
notations above, it is $\tilde \tau (L_V)$. 
Consider the relative tangent bundle with the natural morphism:
$$(\tau f_V)^\ast N_{\tilde \tau  (L_V) / L 
\times_{\mathbb P^2} X } 
= f_V^\ast (\tau^\ast T_{X / \mathbb P^2}) 
\rightarrow
(\sigma, \sigma \circ i)^\ast T_{{\mathcal E} \times i^\ast {\mathcal E} / 
\bar V} 
\vert_{L_V} $$


Since $L_V \times_{\mathbb P^2} X$ is the blowup of ${\mathcal E}_{L_V}
\times_{L_V} i^\ast {\mathcal E}_{L_V}$, by the behavior of the 
canonical class under blow-ups, we get the formula
\begin{equation} \label{e1}
(\sigma, \sigma \circ i)^\ast 
\omega_{{\mathcal E} \times_{\bar V} i^\ast {\mathcal E} / \bar V} \vert_{L_V} \cong 
\tilde\tau^\ast \omega_{X / \mathbb P^2} (-L_V \cap \tilde D) \cong 
f^\ast [ {\mathcal O}_{\mathbb P^2} (3) (-\frac{1}{2} D)] \vert_{L_V},
\end{equation} 
and the last isomorphism holds because $f$ is 2:1.

\medskip

However, 
$\omega_{{\mathcal E} \times_{\bar V} i^\ast {\mathcal E} / \bar V} = \pi_1^\ast 
\omega_{{\mathcal E} / \bar V} \otimes \pi_2^\ast \omega_{i^\ast {\mathcal E} / \bar V}$, 
where $\pi_1$ and $\pi_2$ are the natural projections. 
Denote by $\bar {\mathcal M}_{1,1}$ the coarse moduli space of marked elliptic 
curves. The singular locus of a family of elliptic curves maps to the 
boundary of $\bar {\mathcal M}_{1,1}$. If we denote the pull-back of the 
boundary of $\bar {\mathcal M}_{1,1}$ by $\delta$, 
then it is well known  that 
$\omega_{{\mathcal E} / \bar V} = \frac{\delta}{12}$  
(see \cite[Section 3.E]{hm}).


\medskip

Therefore, using the classifying map from $L_V$ to the moduli space 
$\bar {\mathcal M}_{1,1}$, we obtain 
\begin{equation} \label{e2}
(\sigma, \sigma \circ i)^\ast 
\omega_{{\mathcal E} \times i^\ast {\mathcal E} / \bar V} \cong 
{\mathcal O}_{\bar V} \Big ( \frac{G + i^{-1}G}{12} \Big ) = 
{\mathcal O}_{\bar V} \Big ( \frac{f^{-1}(f(G))}{12} \Big ) = 
f^\ast {\mathcal O}_{\mathbb P^2} \Big ( \frac{f(G)}{12} \Big ) 
\end{equation}
which is well-defined on $L_V$. 

When we compare the isomorphisms (\ref{e1}) and (\ref{e2}), we get: 
$$f^\ast {\mathcal O}_{\mathbb P^2} \Big ( \frac{f(G)}{12} \Big ) \vert_{L_V} 
\cong 
f^\ast [ {\mathcal O}_{\mathbb P^2} (3) (-\frac{1}{2} D)] \vert_{L_V},$$
or equivalently,
$$f^\ast {\mathcal O}_{\mathbb P^2} \Big ( \frac{D}{2} + \frac{f(G)}{12} \Big )
\vert_{L_V} \cong 
f^\ast {\mathcal O}_{\mathbb P^2} (3) \vert_{L_V}$$

Comparing the degrees, we obtain the relation: 
$$\frac{1}{2} \text{deg} (D) + \frac{1}{12} \text{deg} f(G) = 3$$
The degrees of $D$ and $f(G)$ are positive even integers, 
hence there are two possibilities: 
$(\text{deg}(D), \text{deg}(G)) = (2,24)$ or $(4,12)$.

\medskip

{\bf Case 1:} $(\text{deg}(D), \text{deg}(f(G))) = (4,12)$

\medskip

First we consider the case when $D$ is smooth. 
Since $\text{deg}(D)=4$, $\bar V$ is a del Pezzo surface
($K_{\bar V} < 0$). We want to show that 
${\mathcal E}$ is rationally connected. Take two general points $p, q \in 
{\mathcal E}$. Then $f(\pi(p)), f(\pi(q))$ are two general points in 
$\mathbb P^2$. Then $f(\pi(p)) \in L_p$, where $L_p$ is a tangent line to 
$D$ passing through $f(\pi(p))$, and $f(\pi(q)) \in L_q$, where $L_q$ is 
also a tangent line. Let $L_p \cap L_q = \{ r \}$.

\medskip

Take a tangent line $L$ to $D$ and pull it back to $\bar V$: 
$L_{\bar V} \doteq L \times_{\mathbb P^2} \bar V$. Its normalization 
is $\tilde L_{\bar V}$ and ${\mathcal E} \vert_{\tilde L_{\bar V}}$ is an 
elliptic fibration over $\mathbb P^1$ with 12 nodal fibers. 
The surface is rational, because it is deformation equivalent to $\mathbb P^2$ 
blown-up at 9 points in the base locus of a pencil of plane cubics 
(\cite[Section 5.12]{bpv}).

\medskip

We can lift the lines $L_p$ and $L_q$ 
to ${\mathcal E}$ and get ${\mathcal E} \vert_{\tilde L_{p, \bar V} }$ and 
 ${\mathcal E} \vert_{\tilde L_{q, \bar V} }$ which are rational surfaces. We 
can connect any two points on a rational surface with a rational curve. 
Connect $p$ to $\tilde r$ and $q$ to $\tilde r$, where 
$\tilde r \in (f\pi)^{-1} (r)$, so $p$ is rationally chain connected to $q$.

\medskip

The point $\tilde r$ is a smooth point of ${\mathcal E}$. In characteristic 0, 
we can smooth the nodal rational curve if we have a general pair $p, q$ and 
if $\tilde r$ is smooth (\cite[Section 4.6]{de}). 
Therefore, ${\mathcal E}$ is rationally connected.

\medskip

Fix a section $s$ of $i^\ast {\mathcal E} \dashrightarrow \bar V$. Then we get a 
section $\tilde s$ of ${\mathcal E} \times_{\bar V} i^\ast {\mathcal E} 
\dashrightarrow {\mathcal E}$. Since we have a finite morphism 
${\mathcal E} \times_{\bar V} i^\ast {\mathcal E} \dashrightarrow X$, we get a 
finite morphism ${\mathcal E} \dashrightarrow X$.

\medskip

The image of a finite morphism from a rationally connected variety 
to a hyper-K\"ahler manifold is of dimension at most 
$\frac{1}{2} \text{dim}(X)$, because we have a $(2,0)$ form on $X$ and 
there are no holomorphic $(2,0)$ or $(1,0)$ forms on a rational variety. 
However, $\text{dim}({\mathcal E}) = 3$ and 
it is bigger than $\frac{1}{2} \text{dim}(X) = 2$ - a contradiction.

\medskip

Now assume $D$ is singular. Consider the linear system of lines $l$ 
containing a fixed singular point $r$ of $D$. Consider the associated 
linear system $f^{-1}(l)$. These divisors will typically be all singular 
because $\bar V$ is singular at $f^{-1}(r)$. However, the linear system 
of strict transforms $\widetilde{f^{-1}(l)}$ on $Bl_{f^{-1}(r)}{\bar V}$ 
is a basepoint free pencil of divisors on the normal surface 
$Bl_{f^{-1}(r)}{\bar V}$. Thus, by Bertini's theorem, a general member 
of this pencil is smooth and intersects $G$ transversally. In particular, 
the surface ${\mathcal E} \vert_{\widetilde{f^{-1}(l)}}$ is smooth for 
such a member.

\medskip

Let $L$ be any line passing through a singular point $r \in D$. 
Every component of $f^{-1}(L)$ is a rational curve. Without loss of 
generality we shall consider the irreducible case. The normalization 
$\widetilde{f^{-1}(L)}$ is isomorphic to $\mathbb P^1$. Therefore, 
${\mathcal E} \vert_{\widetilde{f^{-1}(L)}}$ is an elliptic surface 
over $\mathbb P^1$. Since there is the following relation of 
intersection numbers: $f^{-1}(L) \cdot G = L \cdot f(G) = 12$, 
it follows that ${\mathcal E} \vert_{\widetilde{f^{-1}(L)}}$ is a 
rational surface.

\medskip

Consider the 1-parameter family $\mathcal F$ of such rational surfaces 
parametrized by $\mathbb P^1$ (since the line $L$ varies in $\mathbb P^1$). 
The 3-fold $\mathcal F$ is rationally connected inside a 4-dimensional 
hyper-K\"ahler manifold which is impossible, because 
$\text{dim}(\mathcal F) > \frac{1}{2} \text{dim}(X)$ 
and there are no holomorphic $(2,0)$ forms on $\mathcal F$.

\medskip

We ruled out the first case completely and the only remaining case is:

\medskip

{\bf Case 2:} $(\text{deg}(D), \text{deg}(f(G))) = (2,24)$

\medskip

First we consider the case when $D$ is smooth. 
Then $\bar V \cong \mathbb P^1 \times \mathbb P^1$.

\medskip

Take a tangent line $L$ to the conic and pull it back to $\bar V$: 
$L_{\bar V} \doteq L \times_{\mathbb P^2} \bar V$. Its normalization 
$\tilde L_{\bar V}$ is reducible and consists of two copies of $\mathbb P^1$, 
say $\tilde L_{\bar V, 1}$ and $\tilde L_{\bar V, 2}$.

\medskip

{\bf Case 2.1:} \\
${\mathcal E} \vert_{\tilde L_{\bar V, k} }$ is an elliptic fibration over 
$\mathbb P^1$ with 12 nodal fibers, $k = 1,2$. 
Then we repeat the same argument as in Case 1 in order to exclude this case.

\medskip

{\bf Case 2.2:} \\
${\mathcal E} \vert_{\tilde L_{\bar V, 1} }$ is an elliptic fibration over 
$\mathbb P^1$ with 24 nodal fibers and ${\mathcal E} 
\vert_{\tilde L_{\bar V, 2}}$ 
is an elliptic fibration over $\mathbb P^1$ with no singular fibers. 
Then  ${\mathcal E} \vert_{\tilde L_{\bar V, 2}}$ is the trivial fibration and 
therefore, ${\mathcal E}$ is the pull back of an elliptic fibration on 
$\mathbb P^1$ through the projection on this factor. And, since the elliptic 
fibration on $\mathbb P^1$ has 24 nodal fibers, it is an elliptic K3 surface 
$S \rightarrow \mathbb P^1$.

\medskip

After considering all the cases, we see that: 
$${\mathcal E} \times_{(\mathbb P^1 \times \mathbb P^1)} i^\ast {\mathcal E} 
 = (S \times \mathbb P^1) \times_{(\mathbb P^1 \times \mathbb P^1)}
( \mathbb P^1 \times S).$$
We want to prove that 
$(S \times \mathbb P^1) \times_{(\mathbb P^1 \times \mathbb P^1)}
( \mathbb P^1 \times S) \cong S \times S$. Indeed, we have the following 
commutative fiber diagram:

$$\begin{array}{ccc}
S \times S & \longrightarrow & \mathbb P^1 \times S \\ 
\downarrow  &  \searrow    & \downarrow \\
S \times \mathbb P^1 & \longrightarrow & \mathbb P^1 \times \mathbb P^1 \\
\end{array}$$

Therefore, 
$${\mathcal E} \times_{(\mathbb P^1 \times \mathbb P^1)} i^\ast {\mathcal E} 
 = (S \times \mathbb P^1) \times_{(\mathbb P^1 \times \mathbb P^1)}
( \mathbb P^1 \times S) \cong S \times S.$$

By construction, $X$ is birational to the desingularization of 
${\mathcal E} \times_{(\mathbb P^1 \times \mathbb P^1)} i^\ast {\mathcal E} / 
{\tilde i}$  which is birational to 
$S \times S / {\mathbb Z_2}$. Therefore, $X$ is birational 
to the desingularization of $S \times S / {\mathbb Z_2}$ which is 
$\text{Hilb}^2(S)$
by Fogarty's theorem (Theorem \ref{fog}).

\medskip

Now consider the case when $D$ is singular. Since $\bar V$ is normal, 
$D$ cannot be a double line. Therefore, $D$ is the union of two lines and 
$\bar V$ is a singular quadric cone.

\medskip

Let $L$ be a line in ${\mathbb P}^2$ passing through the node 
$r \in D$. Then $f^{-1}(L)$ is the union of two lines $L_1$ and $L_2$ 
each one of which is a line in the cone $\bar V$ passing through its 
vertex. On a quadric cone all lines are algebracially equivalent, and in 
particular $L_1$ and $L_2$ are algebraically equivalent. Since 
$\text{deg}(G)=24$, $L_1 \cdot G = L_2 \cdot G = 12$. 
Therefore, the surface ${\mathcal E} \vert_{L_i}$ is rational (for $i=1,2$). 
We constructed a rationally parametrized 1-parameter family of rational 
surfaces which cannot exist in a 4-dimensional hyper-K\"ahler 
manifold. With this we finish the proof of our main theorem. 
\hfill $\Box$
\end{proof}

\paragraph*{Acknowledgements.}
The author is thankful to Prof. G. Tian for 
asking questions leading to this problem. She thanks the referees for 
their helpful suggestions and for improving the exposition of the paper. 
The author is grateful to the Institute for Advanced Study for their 
hospitality.

\providecommand{\bysame}{\leavevmode\hbox to3em{\hrulefill}\thinspace}
%
%

\bibliographystyle{amsalpha}
\bibliographymark{References}
\def\cprime{$'$}

\end{document}